\documentclass[12pt,twoside]{article}    
\usepackage{amsfonts,amsmath,latexsym}

\newtheorem{thm}{Theorem}[section]
\newtheorem{rk}[thm]{Remark}
\newtheorem{lem}[thm]{Lemma}
\newtheorem{cor}[thm]{Corollary}

\newtheorem{ex}[thm]{Example}
\numberwithin{equation}{section}
\newenvironment{proof}[1][Proof]{\textbf{#1.} }{\ \rule{0.5em}{0.5em}}

\makeatletter \@addtoreset{equation}{section} \makeatother
\begin{document}


\pagestyle{myheadings}

\markboth{\hfill {\small AbdulRahman Al-Hussein } \hfill}{\hfill
{\small Necessary and sufficient conditions of optimal control for SDEs} \hfill }


\thispagestyle{plain}


\begin{center}
{\large \textbf{Necessary and sufficient conditions of optimal control for infinite dimensional SDEs}$^{*}$\footnotetext{$^{*}$
This work is supported by King Abdulaziz City for Science and Technology (KACST), Riyadh, Saudi Arabia. }} \\
\vspace{0.7cm} {\large AbdulRahman Al-Hussein }
\\
\vspace{0.2cm} {\footnotesize
{\it Department of Mathematics, College of Science, Qassim University, \\
 P.O.Box 6644, Buraydah 51452, Saudi Arabia \\ {\emph E-mail:} hsien@qu.edu.sa, alhusseinqu@hotmail.com}}
\end{center}

\begin{abstract}
A general maximum principle (necessary and sufficient conditions) for an optimal control problem governed by a stochastic differential equation driven by an infinite dimensional martingale is established. The solution of this equation takes its values in a separable Hilbert space and the control domain need not be convex. The result is obtained by using the adjoint backward stochastic differential equation.
\end{abstract}



\section{Introduction}\label{sec1}
This paper studies the following form of a controlled stochastic differential equation (SDE in short):
\begin{eqnarray}\label{intr:sde}
 \left\{ \begin{array}{ll}
              d X(t) = F (X(t) , u (t)) d t + G (X (t)) d M(t) , \;\; 0 \leq t \leq T,\\
             \; X(0) = x_0 ,
         \end{array}
 \right.
\end{eqnarray}
where $M$ is a continuous martingale taking its values in a separable Hilbert space $ K ,$ while $F ,$ $G$ are some mappings with properties to be given later and $u (\cdot )$ represents a control variable. We will be interested in minimizing the cost functional:
\[
J(u (\cdot ) ) = \mathbb{E} \, [ \, \int_0^T  \ell ( X^{u (\cdot ) } (t) , u (t) ) \, dt + h ( X^{u (\cdot ) } (T) )
\, ] \]
over a set of admissible controls.

We shall follow mainly the ideas of Bensoussan in \cite{[Be82]}, \cite{[Be-book]}, Zhou in \cite{[Zhou93]}, \cite{[Y-Z]}, {\O}ksendal et al. \cite{[Oks-Zh05]}, and our earlier work \cite{[Alh-COSA]}. The reader can see our main results in Theorems~(\ref{main thm}, \ref{mainthm2}).

We recall that SDEs driven by martingales are studied in \cite{[Gyo-Kr80]}, \cite{[Kry-Roz07]}, \cite{[Tud89]}, \cite{[Gre-Tu95]} and \cite{[Alh-AMO-10]}. In fact in \cite{[Alh-AMO-10]} we derived the maximum principle (necessary conditions) for optimality of stochastic systems governed by SPDEs. However, the results there show the maximum principle in its local form and also the control domain is assumed to be convex. In this paper we shall try to avoid such conditions as we shall shortly talk about it. Due to the fact that we are dealing here with a non-convex domain of controls, it is not obvious how one can allow the control variable $u(t)$ to enter in the mapping $G $ in (\ref{intr:sde}) and obtain a result like Lemma~\ref{lem:3rd estimate}. This issue was raised also in \cite{[Be82]}. Nevertheless, in some special cases (see \cite{[Alh-Paper10]}) we can allow $G$ to depend on the control, still overcome this difficulty, and prove the maximum principle. The general result is still open as pointed out in \cite[Remark~6.4]{[Alh-AMO-10]}.

The maximum principle in infinite dimensions started after the work of Pontryagin \cite{[Pont61]}. The reader can find a detailed description of these aspects in Li \& Yong \cite{[Li-Yong95]} and the references therein. An expanded discussion on the history of maximum principle can be found in
\cite[P. 153--156]{[Y-Z]}. On the other hand, the use of (linear) backward stochastic differential equations (BSDEs) for deriving the maximum principle for forward controlled stochastic equations was done by Bismut in \cite{[Bi76]}. In this respect, one can see also the works of Bensoussan in \cite{[Be82]} and \cite{[Be-book]}. In 1990 Pardoux \& Peng, \cite{[Pa-Pe90]}, initiated the theory of nonlinear BSDEs, and then Peng studied the stochastic maximum principle in \cite{[Pe90]} and \cite{[Pe-93]}. Since then several works appeared consequently on the maximum principle and its relationship with BSDEs. For example one can see \cite{[H-Pe90]}, \cite{[H-Pe91]}, \cite{[H-Pe96]}, \cite{[Tang-Li]} and \cite{[Y-Z]} and the references of Zhou cited therein. Our earlier work in \cite{[Alh-Stoc09]} has now opened the way to study BSDEs and backward SPDEs that are driven by martingales. One can see \cite{[Marie09]} for financial applications of BSDEs driven by martingales, and \cite{[Bal-Pard05]}, \cite{[Imk10]}, \cite{[Fuh-Tess]} and \cite{[Alh-ROSE11]} for other applications.

In this paper the convexity assumption on the control domain is not required, as we shall consider a suitable perturbation of an optimal control by means of the spike variation method in order to derive the maximum principle in its global form. Then we shall provide sufficient conditions for optimality of our control problem. The results will be achieved mainly by using the adjoint equation of (\ref{intr:sde}), which is a BSDE driven by the martingale $M .$ This can be seen from equation (\ref{adjoint-bse}) in Section~\ref{sec5}. It is quite important to realize that the adjoint equations of such SDEs are in general BSDEs driven by martingales. This happens also even if the martingale $M ,$ which is appearing in equation (\ref{intr:sde}), is a Brownian motion with respect to a right continuous filtration being larger than its natural filtration. There is a discussion on this issue in Bensoussan's lecture note~\cite[Section~4.4]{[Be82]}, and in \cite{[Alh-ROSE06]} and its erratum \cite{[Alh-ROSE10]}.

The paper is organized as follows. Section~\ref{sec2} is devoted to some preliminary notation.
In Section~\ref{sec3} we present our main stochastic control problems.
Then in Section~\ref{sec4} we establish many of our necessary estimates, which will be needed to derive the maximum principle for the control problem of (\ref{intr:sde}). The maximum principle in the sense of Pontryagin for the above control problem is derived in Section~\ref{sec5}. In Section~\ref{sec6} we establish some sufficient conditions for optimality for this control problem, and present some examples as well.

\section{Preliminary notation}\label{sec2}
Let $(\Omega , \mathcal{F} , \mathbb{P} )$ be a complete probability
space, filtered by a continuous filtration $\{\mathcal{F}_t
\}_{t \geq 0 } ,$ in the sense that every square integrable $K$-valued
martingale with respect to $\{\mathcal{F}_t \, , \; 0 \leq t \leq T \}$ has a continuous version.

Denoting by $\mathcal{P}$ the predictable $\sigma$\,-\,algebra of subsets of $\Omega\times
[0 , T]$ we say that a $K$\,-\,valued process is \emph{predictable} if it is
$\mathcal{P}/\mathcal{B}(K)$ measurable. Suppose that $\mathcal{M}^2_{[0 ,
T]} (K)$ is the Hilbert space of cadlag square integrable
martingales $\{ M (t) , 0 \leq t \leq T \} ,$ which take their values in
$K .$ Let $\mathcal{M}^{2 , c}_{[0 , T ]} (K) $ be the subspace of
$\mathcal{M}^2_{[0 , T]} (K)$ consisting of all continuous square
integrable martingales in $K .$ Two
elements $M$ and $N$ of $\mathcal{M}^2_{[0 , T]} (K) $ are said to be {\it very strongly orthogonal} (or shortly VSO) if $$\mathbb{E}\, [ M
(\tau )\otimes N (\tau ) ] =\, \mathbb{E}\, [ M (0)\otimes N (0) ] ,$$ for
all $[ 0 , T]$\,-\,valued stopping times $\tau .$

Now for $M \in \mathcal{M}^{2 , c}_{[0 , T ]} (K) $ we shall use the notation $<M>$ to mean
the predictable quadratic variation of $M $ and
similarly $<<M>>$ means the predictable tensor quadratic variation of $M ,$ which takes its values
in the space $L_1(K)$ of all nuclear operators on $K .$ Precisely, $M\otimes M - <<M>> \, \in \mathcal{M}^{2 , c}_{[0 , T ]} (L_1 (K)) .$
We shall assume for a given fixed $M \in \mathcal{M}^{2 , c}_{[0 , T ]}(K)$ that there exists a measurable mapping ${\mathcal{Q} (\cdot) : [ 0 , T ] \times \Omega \rightarrow L_1 (K)}$ such that $\mathcal{Q} (t)$ is symmetric,
positive definite, $\mathcal{Q} (t) \leq \mathcal{Q} $ for some positive definite
nuclear operator $\mathcal{Q}$ on $K ,$ and satisfies the following equality:
$$<< M >>_t \; = \int_0^t \mathcal{Q} (s) \, ds .$$ We refer the reader to Example~\ref{Example1} for a precise computation of this process $\mathcal{Q} (\cdot ) . $

For fixed $(t , \omega ), $ we denote by
$L_{\mathcal{Q}(t , \omega )} (K)$ to the set of all linear operators $\varphi : \mathcal{Q}^{1/2}(t ,
\omega ) (K) \rightarrow K $ and satisfy $\varphi \mathcal{Q}^{1/2}(t , \omega)
\in L_2 (K) ,$ where $L_2 (K)$ is the space of all Hilbert-Schmidt operators from
$K$ into itself. The inner product and norm in $L_2 (K)$  will be denoted respectively by $\big{<} \cdot , \cdot \big{>}_2 $
and $|| \cdot ||_2$. Then the stochastic integral $\int_0^{\cdot} \Phi (s) d M(s) $ is defined
for mappings $\Phi$ such that for each $(t , \omega ) , \; \Phi (t ,
\omega ) \in L_{\mathcal{Q} (t , \omega )} (K) , \; \Phi \mathcal{Q}^{1/2} (t, \omega) (h) \; \forall \; h \in K$ is predictable, and
\[ \mathbb{E}\; [ \, \int_0^T || ( \Phi \mathcal{Q}^{1/2} ) (t) ||_2^2
\; dt \, ] < \infty . \]

Such integrands form a Hilbert space with respect to the
scalar product $( \Phi_1 , \Phi_2 ) \mapsto \mathbb{E}\; [ \,
\int_0^T \big{<} \Phi_1 \mathcal{Q}^{1/2} (t)\; ,
\Phi_2 \mathcal{Q}^{1/2} (t)\big{>} \; dt \, ] .$ Simple processes taking values in $L(K ; K)$
are examples of such integrands. By letting $\Lambda^2 ( K ; \mathcal{P} , M ) $ be the
closure of the set of simple processes in this Hilbert space, it becomes a Hilbert subspace.
We have also the following isometry property:
\begin{eqnarray}\label{isometry property}
\mathbb{E} \, \Big[ \, | \int_0^T \Phi (t) dM(t) |^2 \Big] &=&
\mathbb{E} \, \Big[ \, \int_0^T || \Phi (t) \mathcal{Q}^{1/2} (t)
||_2^2 \, ds \Big]
\end{eqnarray}
for mappings $\Phi \in \Lambda^2 ( K ; \mathcal{P} , M ) .$ For more details and proofs we refer the reader to \cite{[Me-P]}.

On the other hand, we emphasize that the process $\mathcal{Q} (\cdot)$ will be play an important role in deriving the adjoint equation of the SDE~(\ref{intr:sde}) as it can be seen from equations (\ref{defn of H}), (\ref{adjoint-bse}) in Section~\ref{sec5}. This is due to the fact that the integrand $\Phi$ is not necessarily bounded. More precisely, it is needed in order for the mapping $\nabla_{x} H ,$ which appear in both equations, to be defined on the space $L_2 (K) ,$ since the process $Z^{u (\cdot )}$ there need not be bounded. This always has to be considered when working with BSDEs or BSPDEs driven by infinite dimensional martingales.

Next let us introduce the following space:
\[ L^2_{\mathcal{F}} ( 0 , T ;  E  ) := \{ \psi : [ 0 , T]\times \Omega \rightarrow E,
\; \text{predictable and} \; \mathbb{E} \, [ \int_0^T | \psi (t) |^2 d
t \, ] < \infty \, \} ,\] where $E$ is a separable Hilbert
space.

Since $\mathcal{Q} (t) \leq \mathcal{Q} $ for all $t \in [ 0 , T]$ a.s., it follows from \cite[Proposition~2.2]{[Alh-Int10]} that if $\Phi \in L^2_{\mathcal{F}} ( 0 , T ; L_{\mathcal{Q}} (K) )  $ (where as above $L_{\mathcal{Q}} (K) ) = L_2 (\mathcal{Q}^{1/2} (K) ; K) ) ,$ the space of all Hilbert-Schmidt operators from $\mathcal{Q}^{1/2} (K)$ into $K $), then $\Phi \in \Lambda^2 ( K ; \mathcal{P} , M ) $ and
\begin{eqnarray}\label{eq:boundedness of Q}
\mathbb{E}\; [ \, \int_0^T || \Phi (t) \mathcal{Q}^{1/2} (t) ||^2_2 \; dt \, ] \leq \mathbb{E}\; [ \, \int_0^T || \Phi (t) ||_{L_{\mathcal{Q}}(K)}^2 \; d t \, ] .
\end{eqnarray}
An example of such a mapping $\Phi$ is the mapping $G$ in equation (\ref{intr:sde}); see the domain of $G$ in the introduction of the following section.

\section{Formulation of the control problem}\label{sec3}
Let $\mathcal{O}$ be a separable Hilbert space and $U$ be a nonempty subset of $\mathcal{O} .$ We say that $u (\cdot ) : [0 , T]\times \Omega \rightarrow \mathcal{O}$ is \emph{admissible} if $u (\cdot ) \in  L^2_{\mathcal{F}} ( 0 , T ; \mathcal{O} )$ and $u (t) \in U \; \; a.e., \; a.s.$
The set of admissible controls will be denoted by $\mathcal{U}_{ad} .$

Let $F: K \times \mathcal{O} \rightarrow K ,$ $G : K \rightarrow L_{\mathcal{Q}} (K) ,$
$\ell: K \times \mathcal{O} \rightarrow \mathbb{R}$ and $h : K \rightarrow \mathbb{R} $ be measurable mappings.
Consider the following SDE:
\begin{eqnarray}\label{forward-see}
 \left\{ \begin{array}{ll}
               d X(t) = F (X(t) , u (t) ) \, d t
              + \, G (X(t) )\, d M(t) , \;\; t \in [0 , T] ,\\
             \; X(0) = x_0 \in K .
         \end{array}
 \right.
\end{eqnarray}
If assumption (E1), which is stated below, holds, then (\ref{forward-see}) attains a unique solution in $L^2_{\mathcal{F}} ( 0 , T ; K ) .$ The proof of this fact can be gleaned from \cite{[Re-Y]} or \cite{[Roz90]}. In this case we shall denote the solution of (\ref{forward-see}) by $X^{u (\cdot ) } .$

\bigskip

Our assumptions are the following. \\
(E1) $F , G , \ell , h$ are continuously Fr\'echet differentiable with respect to
$x ,$ $F$ and $\ell$ are continuously Fr\'echet differentiable with respect to
$u ,$ the derivatives $F_x , \, F_{u } , \,  G_x , \ell_x , \ell_{u} $ are uniformly bounded, and
$$ | h_x |_{L(K; K)} \leq k \, ( 1 + |x|_{K} ) $$ for some constant $k > 0 .$

In particular, $|F_x|_{L(K,K)} \leq C_1 , \, ||G_x||_{L(K,L_{\mathcal{Q}}(K))}  \leq C_2 , \,  |F_{v}|_{L(\mathcal{O},K)} \leq C_3 ,  $ for some positive constants $C_i , \; i =1, 2, 3 ,$ and similarly for $\ell .$

\noindent (E2) $\ell_x$ satisfies Lipschitz condition with respect to $u $ uniformly in $x .$

\bigskip

Consider now the \emph{cost functional}:
\begin{equation}\label{cost functional}
J(u (\cdot ) ) : = \mathbb{E} \, [ \, \int_0^T  \ell ( X^{u
(\cdot ) } (t) , u (t) ) \, dt + h ( X^{u (\cdot ) } (T) )
\, ] ,
\end{equation}
for $u (\cdot ) \in \mathcal{U}_{ad} .$

\smallskip

The control problem here is to minimize (\ref{cost functional}) over the set $\mathcal{U}_{ad} .$
Any $u^{*} ( \cdot ) \in \mathcal{U}_{ad} $ satisfying
\begin{equation}\label{value-function}
J(u^{*} ( \cdot ) ) = \inf \{ J(u (\cdot ) ): \; u (\cdot ) \in \,
\mathcal{U}_{ad} \}
\end{equation}
is called an \emph{optimal control}, and its corresponding solution $X^* := X^{u^{*} (\cdot ) }$ to (\ref{forward-see}) is called an \emph{optimal solution} of the stochastic optimal control problem (\ref{forward-see})-(\ref{value-function}). In this case the pair $( X^{*} \, , u^{*} (\cdot ) )$ in this case is called an \emph{optimal pair}.

\begin{rk}\label{rk:dependence on t}
We mention here that the mappings $F , G$ and $\ell$ in (\ref{forward-see}) and (\ref{cost functional}) can be taken easily to depend on time $t$ with a similar proof as established in the following sections, but rather, having more technical computations.
\end{rk}

Since this control problem has no constraints we shall deal generally with progressively measurable controls. However, for the case when there
are final state constraints, one can mimic our results in Sections (\ref{sec4}, \ref{sec5}, \ref{sec6}), and use Ekeland's variational principle in a similar way to \cite{[Maz]}, \cite{[Pe90]} or \cite{[Y-Z]}.

In the following section we shall begin with some variational method in order to derive our main variational inequalities that are necessary to establish the main result of Section~\ref{sec5}.

\section{Estimates}\label{sec4}
Let $(X^{*} , u^{*} (\cdot ) )$ be the given optimal pair. Let $0 \leq t_0 < T$ be fixed and $0 < \varepsilon < T - t_0 .$ Let $v$ be a random variable taking its values in $U ,$  $\mathcal{F}_{t_0}$\,-\,measurable and $\displaystyle{\sup_{\omega \in \Omega }} \, | v (\omega ) | < \infty .$ Consider the following spike variation of the control $u^{*} (\cdot )$:
\begin{eqnarray}\label{eq:spike variation}
 u_{\varepsilon } (t) = \left\{ \begin{array}{ll}
u^{*} (t) \;\; & \text{ if} \; \; t \in [0 , T] \backslash [t_0 , t_0 + \varepsilon ] \\
v \;\; & \text{ if} \; \; t \in [t_0 , t_0 + \varepsilon ] .
\end{array}
 \right.
\end{eqnarray}

Let $X^{u_{\varepsilon } (\cdot ) }$ denote the solution of the SDE~(\ref{forward-see}) corresponding to $ u_{\varepsilon } (\cdot ) .$ We shall denote it briefly by $X_{\varepsilon } .$ Observe that $X_{\varepsilon } (t) = X^{*} (t) $ for all $0 \leq t \leq t_0 .$

The following lemmas will be very useful in proving the main results of Section~\ref{sec5}.
\begin{lem}\label{lem:1st estimate}
Let (E1) hold. Assume that $\{ p(t) , \; t_0 \leq t \leq T \}$ is the solution of the following linear equation:
\begin{eqnarray}\label{eq:p}
 \left\{ \begin{array}{ll}
               d p(t) = F_x (X^{*}(t) , u^{*} (t) ) \, p (t) \, dt + G_x (X^{*} (t) )\, p (t) \, d M(t) , \hspace{0.5cm} t_0 < t \leq T,\\
               \; p (t_0 ) = F (X^{*}(t_0) , v ) - F (X^{*}(t_0) , u^{*} (t_0) ) .
         \end{array}
 \right.
\end{eqnarray}
Then
\[
\sup_{t \in [t_0 , T]} \mathbb{E} \, [ \, | p (t) |^2 \, ] \, < C
\]
for some positive constant $C .$
\end{lem}
\begin{proof}
With the help of (E1) apply It\^{o}'s formula to compute $| p (t) |^2 ,$ and take the expectation. The required result follows then by using Gronwall's inequality.
\end{proof}

\begin{lem}\label{lem:2nd estimate}
Assuming (E1) we have
\[
\mathbb{E} \, [ \, \sup_{t_0 \leq t \leq T } \, | X_{\varepsilon} (t) - X^{*} (t) |^2 \, ] = o (\varepsilon) .
\]
\end{lem}
\begin{proof}
For $t_0 \leq t \leq t_0 + \varepsilon $ one observes that
\begin{eqnarray}\label{eq:est1}
X_{\varepsilon} (t) - X^{*} (t) &=& \int_{t_0}^t [ F( X_{\varepsilon}(s) , v ) - F( X^{*}(s) , v ) ] \, ds \nonumber \\ && + \, \int_{t_0}^t [ F( X^{*}(s) , v ) - F( X^{*}(s) , u^{*} (s) ) ] \, ds \nonumber \\ &&
+ \, \int_{t_0}^t [ G( X_{\varepsilon}(s) ) - G( X^{*}(s) ) ] dM(s),
\end{eqnarray}
or, in particular,
\begin{eqnarray}\label{eq:est1a}
&& \hspace{-1cm} | X_{\varepsilon} (t) - X^{*} (t) |^2 \leq 3\, (t-t_0 ) \int_{t_0}^t | F( X_{\varepsilon}(s) , v ) - F( X^{*}(s) , v ) |^2 \, ds \nonumber \\ &&
\hspace{3cm} + \, 3 \, (t-t_0 ) \int_{t_0}^t | F( X^{*}(s) , v ) - F( X^{*}(s) , u^{*} (s) ) |^2 \, ds \nonumber \\ &&
\hspace{4cm} + \, 3 \; | \int_{t_0}^t [ G( X_{\varepsilon}(s) ) - G( X^{*}(s) ) ] dM(s) |^2 .
\end{eqnarray}

But Taylor expansion implies the three identities:
\begin{eqnarray}\label{eq:est1a-F-x}
&& \hspace{-1.5cm} F( X_{\varepsilon}(s) , v ) - F( X^{*}(s) , v ) \nonumber \\ && \hspace{-0.5cm} = \, \int_0^1 F_x ( X^{* }(s) , u^{*} (s) +\lambda ( X_{\varepsilon }(s) - X^{*}(s) )) \, ( X_{\varepsilon }(s) - X^{*}(s) ) \, d\lambda ,
\end{eqnarray}
\begin{eqnarray}\label{eq:est1a-F-u}
&& \hspace{-1cm} F( X^{*}(s) , v ) - F( X^{*}(s) , u^{*} (s) ) \nonumber \\ && \hspace{1cm} = \, \int_0^1 F_v ( X^{* }(s) , u^{*} (s) + \lambda (v - u^{*} (s) )) \, (v - u^{*} (s) ) \, d\lambda ,
\end{eqnarray}
and
\begin{eqnarray}\label{eq:est1a-G-1}
&& \hspace{-1.5cm} G( X_{\varepsilon}(s) ) - G( X^{*}(s) ) \nonumber \\ && \hspace{0.5cm} = \, \int_0^1 G_x ( X^{*} (s) + \lambda ( X_{\varepsilon }(s) - X^{*}(s) )) \, ( X_{\varepsilon }(s) - X^{*}(s) ) \, d\lambda \nonumber \\ && \hspace{3in} =: \Phi (s) \; ( \in L_{\mathcal{Q}} (K) ).
\end{eqnarray}

Then, by using (\ref{eq:est1a-G-1}), the isometry property (\ref{isometry property}), (\ref{eq:boundedness of Q}) and (E1) we deduce that for all $t \in [ t_0 , t_0 + \varepsilon ] ,$
\begin{eqnarray}\label{eq:est1a-G-2}
&& \hspace{-1cm} \mathbb{E}\, [ \, | \int_{t_0}^t \big( G( X_{\varepsilon}(s) ) - G( X^{*}(s) ) \big) dM(s) |^2 \, ] = \mathbb{E}\, [ \, | \int_{t_0}^t \Phi (s) dM(s) |^2 \, ]  \nonumber \\
&& \hspace{-1cm} = \; \mathbb{E}\, [ \, \int_{t_0}^t || \Phi (s) \mathcal{Q}^{1/2} (s) ||_2^2 \, ds \, ] \nonumber \\
&& \hspace{-1cm} \leq \; \mathbb{E}\, [ \, \int_{t_0}^t || \Phi (s) ||_{L_{\mathcal{Q}} (K)}^2 \, ds \, ] \nonumber \\
&& \hspace{-1cm} = \; \mathbb{E}\, [ \, \int_{t_0}^t || \int_0^1 G_x ( X^{*} (s) + \lambda ( X_{\varepsilon }(s) - X^{*}(s) )) \, ( X_{\varepsilon }(s) - X^{*}(s) ) \, d\lambda ||_{L_{\mathcal{Q}} (K)}^2 \, ds \, ] \nonumber \\
&& \hspace{-1cm} \leq \; \mathbb{E}\, [ \, \int_{t_0}^t \int_0^1 || G_x ( X^{*} (s) + \lambda ( X_{\varepsilon }(s) - X^{*}(s) )) \, ( X_{\varepsilon }(s) - X^{*}(s) ) ||_{L_{\mathcal{Q}} (K)}^2  \, d\lambda \, ds \, ] \nonumber \\
&& \hspace{-1cm} \leq \; C_2 \; \mathbb{E}\, [ \, \int_{t_0}^t | X_{\varepsilon }(s) - X^{*}(s) |^2  \, ds \, ] .
\end{eqnarray}

Therefore, from (\ref{eq:est1a}), (\ref{eq:est1a-F-x}), (\ref{eq:est1a-F-u}), (E1) and (\ref{eq:est1a-G-2}), it follows evidently that
\begin{eqnarray*}
&& \hspace{-1cm} \mathbb{E}\, [ \, | X_{\varepsilon} (t) - X^{*} (t) |^2 \, ]
\leq  3\, \big{(} C_1 \, (t-t_0 ) + C_2 \big{)} \int_{t_0}^t \mathbb{E}\, [ \, |\, X_{\varepsilon} (s) - X^{*} (s) \, |^2  \, ] \, ds
\nonumber \\ && \hspace{2in} + \, 3\, (t-t_0 ) \, C_3 \int_{t_0}^t \mathbb{E}\, [ \, |\, v - u^{*} (s) \, |^2 \, ] \, ds ,
\end{eqnarray*}
for all $t \in [ t_0 , t_0 + \varepsilon ] .$

Hence by using Gronwall's inequality we obtain
\begin{eqnarray}\label{eq:est2}
&& \hspace{-2cm}\mathbb{E}\, [ \, |\, X_{\varepsilon} (t) - X^{*} (t) \, |^2 \, ]  \leq 3\, C_3 \, (t-t_0 ) \, e^{3\, \big{(} C_1 \, (t-t_0 ) + C_2 \big{)} (t-t_0 ) } \, \nonumber \\ &&
\hspace{2in} \times \int_{t_0}^{t_0 + \varepsilon } \mathbb{E}\,  [ \, | v - u^{*} (s) |^2 \, ] \, ds ,
\end{eqnarray}
for all $t \in [t_0 , t_0 + \varepsilon ] .$ Consequently,
\begin{eqnarray}\label{eq:est3}
&& \hspace{-2cm} \mathbb{E}\, [ \, \int_{t_0}^{t_0 + \varepsilon } |\, X_{\varepsilon} (t) - X^{*} (t) \, |^2 \, dt \, ]  \leq 3\, C_3 \,  \varepsilon^2 \, e^{3 \, (C_1 \,  \varepsilon + C_2 ) \varepsilon } \, \nonumber \\ &&
\hspace{2in} \times \int_{t_0}^{t_0 + \varepsilon } \mathbb{E}\,  [ \, | v - u^{*} (s) |^2 \, ] \, ds .
\end{eqnarray}

It follows then from (\ref{eq:est1a}), (\ref{eq:est2}), standard martingale inequalities, (\ref{eq:est1a-G-2}) and (\ref{eq:est3}) that
\begin{eqnarray}\label{eq:est4}
&& \hspace{-1.5cm} \mathbb{E}\, [ \, \sup_{t_0 \leq t \leq t_0 + \varepsilon } | X_{\varepsilon} (t) - X^{*} (t) |^2 \, ] \nonumber \\ && \hspace{-1.5cm} \leq 3 \, C_3 \, [ \, 3 \, ( C_1 \, \varepsilon + 4 C_2) \, \varepsilon \, e^{3 \, (C_1 \,  \varepsilon + C_2 ) \varepsilon } \, +  \, 1 \, ] \, \varepsilon \, \int_{t_0}^{t_0 + \varepsilon } \mathbb{E}\,  [ \, | v - u^{*} (s) |^2 \, ] \, ds .
\end{eqnarray}

Next, for $t_0 + \varepsilon \leq t \leq T ,$ we have
\begin{eqnarray}\label{eq:est5}
&& X_{\varepsilon} (t) - X^{*} (t) = X_{\varepsilon} (t_0 + \varepsilon ) - X^{*} (t_0 + \varepsilon) \nonumber \\ && \hspace{1.25in} + \, \int_{t_0 + \varepsilon}^t [ F( X_{\varepsilon}(s) , u^{*} (s) ) - F( X^{*}(s) , u^{*} (s) ) ] \, ds \nonumber \\ && \hspace{1.25in}
+ \, \int_{t_0 + \varepsilon }^t [ G( X_{\varepsilon}(s) ) - G( X^{*}(s) ) ] dM(s) .
\end{eqnarray}
Thus by working as before and applying (\ref{eq:est2}) we derive
\begin{eqnarray*}
&& \mathbb{E}\, [ \, \int_{t_0 + \varepsilon }^T |\, X_{\varepsilon} (t) - X^{*} (t) \, |^2 \, dt \, ]  \leq 9 \, C_3 \,  \varepsilon^2 e^{C_4 (\varepsilon )}  \int_{t_0}^{t_0 + \varepsilon } \mathbb{E}\,  [ \, | v - u^{*} (s) |^2 \, ] \, ds
\end{eqnarray*}
and
\begin{eqnarray}\label{eq:est6}
&& \hspace{-0.65cm} \mathbb{E}\, [ \, \sup_{t_0 + \varepsilon \leq t \leq T } | X_{\varepsilon} (t) - X^{*} (t) |^2 \, ]
\leq 27 \, C_3 \,  \varepsilon \, e^{C_4 (\varepsilon )} \, [ \, 1+ \big{(} ( T - t_0 - \varepsilon ) \, C_1 + 4 \, C_2 \big{)} \, \varepsilon \, ] \nonumber \\ &&
\hspace{2.5in} \times \int_{t_0}^{t_0 + \varepsilon } \mathbb{E}\,  [ \, | v - u^{*} (s) |^2 \, ] \, ds ,
\end{eqnarray}
where $C_4 (\varepsilon ) = [ 3 \, \varepsilon^2 + 3 \, ( T - t_0 - \varepsilon )^2 ] \, C_1 + ( T - t_0 + 2\, \varepsilon ) \, C_2 .$

Now (\ref{eq:est4}) and (\ref{eq:est6}) imply that
\begin{eqnarray*}
&& \mathbb{E}\, [ \, \sup_{t_0 \leq t \leq T } | X_{\varepsilon} (t) - X^{*} (t) |^2 \, ]
\leq \big{(} C_5 (\varepsilon ) + C_6 (\varepsilon ) \big{)} \int_{t_0}^{t_0 + \varepsilon } \mathbb{E}\,  [ \, | v - u^{*} (s) |^2 \, ] \, ds ,
\end{eqnarray*}
with the constants $$C_5 (\varepsilon ) = 3 \, C_3 \, [ \, 3 \, ( C_1 \, \varepsilon + 4 C_2) \, \varepsilon \, e^{3 \, (C_1 \,  \varepsilon + C_2 ) \varepsilon } \, +  \, 1 \, ] \, \varepsilon $$
and
$$C_6 (\varepsilon ) = 27 \, C_3 \,  \varepsilon \, e^{C_4 (\varepsilon )} \, [ \, 1+ \big{(} ( T - t_0 - \varepsilon ) \, C_1 + 4 \, C_2 \big{)} \, \varepsilon \, ] .$$
This completes the proof.
\end{proof}

\bigskip
\begin{rk}
We note that for $a.e. \; s, $
\begin{equation}\label{eq:t0}
\frac{1}{\varepsilon } \, \int_{s}^{s + \varepsilon } \mathbb{E} \, [ \, |\phi( X^{*}(t) , u^{*}(t) ) - \phi( X^{*}(s) , u^{*} (s) ) |^2 \, ] \, dt \rightarrow 0 ,  \; \; \text{as} \;\; \varepsilon \rightarrow 0 ,
\end{equation}
for $\phi = F , \ell .$ Indeed, if for example, $\phi = F , $ then we may argue as in (\ref{eq:est1a-F-u}) to see that
\begin{eqnarray}\label{eq:F-t0}
&& \hspace{-0.55cm} \frac{1}{\varepsilon } \, \int_{s}^{s + \varepsilon } \mathbb{E} \, [ \, | F( X^{*}(t) , u^{*}(t) ) - F( X^{*}(s) , u^{*} (s) ) |^2 \, ] \, dt
\nonumber \\ && \hspace{-0.55cm}
= \, \frac{1}{\varepsilon } \, \int_{s}^{s + \varepsilon } \mathbb{E} \, [ \, | \int_0^1 F_v ( X^{* }(t) , u^{*} (s) + \lambda (u^{*} (t) - u^{*} (s) )) \, (u^{*} (t) - u^{*} (s)) \, d\lambda |^2 \, dt \, ]  \nonumber \\ && \hspace{-0.55cm}
\leq \, \frac{1}{\varepsilon } \, \int_{s}^{s + \varepsilon } \mathbb{E} \, [ \, | u^{*} (t) - u^{*} (s) |^2 \, ] \, dt .
\end{eqnarray}
But since $\int_{0}^{T} \mathbb{E} \, [ \, | u^{*} (t) - u^{*} (s) |^2 \, ] \, dt < \infty $ (for fixed $s$), then, as it is well-known from measure theory (e.g. \cite{[Cohn]}), there exists a subset $O$ of $[0 , T]$ such that $\mathbb{P} (O) = 1 $ and the mapping $ O \ni t \mapsto \mathbb{E} \, [ \, | u^{*} (t) - u^{*} (s) |^2 \, ] $ is continuous. Thus, if $s \in O ,$ this function is continuous in a neighborhood of $s ,$ and so we have
\begin{equation*}
\frac{1}{\varepsilon } \, \int_{s}^{s + \varepsilon } \mathbb{E} \, [ \, | u^{*} (t) - u^{*} (s) |^2 \, ] \, dt \rightarrow 0 ,  \; \; \text{as} \;\; \varepsilon \rightarrow 0 ,
\end{equation*}
which by (\ref{eq:F-t0}) implies (\ref{eq:t0}) for $\phi = F .$
\end{rk}

We will choose $t_0$ such that (\ref{eq:t0}) holds for $\phi = F , \ell .$ This assumption will be considered until the end of Section~\ref{sec5}.

\begin{lem}\label{lem:3rd estimate}
Assume (E1). Let $$\xi_{\varepsilon} (t) = \frac{1}{\varepsilon } \, (X_{\varepsilon} (t) - X^{*} (t)) - p (t) , \hspace{0.25cm} t \in [t_0 , T ] .$$ Then
\[  \lim_{\varepsilon \rightarrow 0 } \mathbb{E} \, [ \, | \xi_{\varepsilon} (T) |^2 \, ] \, = \, 0 . \]
\end{lem}
\begin{proof}
First note that, for $t_0 \leq t \leq t_0 + \varepsilon ,$
\begin{eqnarray*}\label{eq1:lem3}
&& d \xi_{\varepsilon } (t) = \frac{1}{\varepsilon } \, [ \, F( X_{\varepsilon }(t) , v ) - F( X^{*}(t) , u^{*} (t) ) - \varepsilon \, F_x (X^{*}(t) , u^{*} (t) ) \, p(t) \, ] dt \nonumber \\ && \hspace{1.5cm}
+ \; \frac{1}{\varepsilon } \, [ \, G( X_{\varepsilon }(t)) - G( X^{*}(t)) - \varepsilon \, G_x (X^{*}(t)) \, p(t) \, ] dM(t) , \nonumber \\
&& \; \xi_{\varepsilon } (t_0 ) = - \big{(} F( X^{*}(t_0) , v ) - F( X^{*}(t_0) , u^{*} (t_0) ) \big{)} .
\end{eqnarray*}
Thus
\begin{eqnarray*}
\xi_{\varepsilon } (t_0 + \varepsilon) &=& \frac{1}{\varepsilon } \, \int_{t_0 }^{t_0 + \varepsilon} [ \, F( X_{\varepsilon }(s) , v ) - F( X^{*}(s) , v) \, ] \, ds  \nonumber \\ && +\; \frac{1}{\varepsilon } \, \int_{t_0 }^{t_0 + \varepsilon} [ \, F( X^{*}(s) , v ) - F( X^{*}(t_0 ) , v) \, ] \, ds  \nonumber \\ && +\; \frac{1}{\varepsilon } \, \int_{t_0 }^{t_0 + \varepsilon} [ \, F( X^{*}(t_0 ) , u^{*} (t_0 ) ) - F( X^{*}(s) , u^{*}(s) ) \, ] \, ds  \nonumber \\ && +\; \frac{1}{\varepsilon } \, \int_{t_0 }^{t_0 + \varepsilon} [ \, G( X_{\varepsilon }(s)) - G( X^{*}(s)) \, ] dM(s) \nonumber \\ &&
- \int_{t_0 }^{t_0 + \varepsilon} F_x (X^{*}(s) , u^{*} (s) ) p(s) \, ds
\nonumber \\ && - \int_{t_0 }^{t_0 + \varepsilon} G_x (X^{*}(s))  p(s) dM(s) .
\end{eqnarray*}

By using (\ref{isometry property}), (\ref{eq:boundedness of Q}) and (E1) we deduce
\begin{eqnarray}\label{eq2:lem3}
&& \hspace{-1.5cm} \mathbb{E}\, [ \, |\, \xi_{\varepsilon } (t_0 + \varepsilon)  \, |^2 \, ] \leq 6\, C_1 \, \mathbb{E}\, [ \, \sup_{t_0 \leq t \leq t_0 + \varepsilon } | X_{\varepsilon } (t) - X^{*} (t) |^2 \, ] \nonumber \\
&& + \, 6 \, \sup_{t_0 \leq t \leq t_0 + \varepsilon } \mathbb{E}\, [ \, | F( X^{* }(t) , v ) - F( X^{*}(t_0 ) , v ) |^2 \, ]
\nonumber \\ && + \, \frac{6}{\varepsilon } \, \int_{t_0 }^{t_0 + \varepsilon} \mathbb{E}\, [ \, | F( X^{* }(s) , u^{*}(s) ) - F( X^{*}(t_0 ) , u^{*}(t_0)  ) |^2 \, ] \, ds \nonumber \\
&& + \, \frac{6 \, C_2}{\varepsilon } \; \mathbb{E}\, [ \, \sup_{t_0 \leq t \leq t_0 + \varepsilon } | X_{\varepsilon } (t) - X^{*} (t) |^2 \, ] \nonumber \\ && + \, 6 \, (C_1 + C_2 ) \; \mathbb{E}\, [ \, \int_{t_0 }^{t_0 + \varepsilon} | p(s) |^2 \, ds \, ] .
\end{eqnarray}
But from (\ref{eq:est4})
\begin{eqnarray}\label{eq3:lem3}
&& \hspace{-1.75cm} \frac{1}{\varepsilon } \; \mathbb{E}\, [ \, \sup_{t_0 \leq t \leq t_0 + \varepsilon } | X_{\varepsilon } (t) - X^{*} (t) |^2 \, ]  \nonumber \\ && \hspace{-1cm}
\leq 3 \, C_3 \, [ \, 3 \, ( C_1 \, \varepsilon + 4 C_2) \, \varepsilon \, e^{3 \, (C_1 \,  \varepsilon + C_2 ) \varepsilon } \, +  \, 1 \, ] \, \int_{t_0}^{t_0 + \varepsilon } \mathbb{E}\,  [ \, | v - u^{*} (s) |^2 \, ] \, ds
 \nonumber \\ && \hspace{3in} \rightarrow 0 \hspace{0.5cm} \text{as} \; \; \varepsilon \rightarrow 0 .
\end{eqnarray}
Also as in (\ref{eq:est1a-F-x}),  by applying (E1) and (\ref{eq:est2}), one gets
\begin{eqnarray}\label{eq4:lem3}
&& \hspace{-1.5cm} \mathbb{E}\, [ \, | F( X^{* }(t) , v ) - F( X^{*}(t_0 ) , v ) |^2 \, ]
\nonumber \\
&& \hspace{-1.5cm} = \, \mathbb{E}\, [ \, | \int_0^1 F_x ( X^{*} (t_0 ) + \lambda ( X^{* }(t) - X^{*}(t_0 ) ) , v ) ( X^{* }(t) - X^{*}(t_0 ) ) \, d\lambda |^2
\nonumber \\
&& \hspace{-1.5cm} \leq \, C_1 \; \mathbb{E}\, [ \, | X^{* }(t) - X^{*}(t_0 ) |^2 \, ] \nonumber \\ &&
\hspace{-1.5cm} \leq 3 \, C_1 \,  C_3 \, \varepsilon \, e^{3\, \big{(} C_1 \, \varepsilon + C_2 \big{)} \varepsilon } \, \int_{t_0}^{t_0 + \varepsilon } \mathbb{E}\,  [ \, | v - u^{*} (s) |^2 \, ] \, ds \rightarrow 0 \hspace{0.5cm} \text{as} \; \; \varepsilon \rightarrow 0 .
\end{eqnarray}

Thus, by applying Lemma~\ref{lem:2nd estimate}, (\ref{eq4:lem3}), (\ref{eq3:lem3}), (\ref{eq:t0}) and Lemma~\ref{lem:1st estimate} in (\ref{eq2:lem3}), we deduce
\begin{eqnarray}\label{eq5:lem3}
&& \hspace{-1.5cm} \mathbb{E}\, [ \, |\, \xi_{\varepsilon } (t_0 + \varepsilon)  \, |^2 \, ] \rightarrow 0 \hspace{0.5cm} \text{as} \; \; \varepsilon \rightarrow 0 .
\end{eqnarray}

Let us now assume that $t_0 + \varepsilon \leq t \leq T .$ In this case we have
\begin{eqnarray*}\label{eq6:lem3}
&& \hspace{-0.5cm} d \xi_{\varepsilon } (t) = \frac{1}{\varepsilon } \, [ \, F( X_{\varepsilon }(t) , u^{*} (t) ) - F( X^{*}(t) , u^{*} (t) ) - \varepsilon \, F_x (X^{*}(t) , u^{*} (t) ) \, p(t) \, ] dt  \nonumber \\
&& \hspace{1.5cm} + \; \frac{1}{\varepsilon } \, [ \, G( X_{\varepsilon }(t)) - G( X^{*}(t)) - \varepsilon \, G_x (X^{*}(t)) \, p(t) \, ] dM(t) ,
\end{eqnarray*}
or, in particular, by setting $$\tilde{\Phi}_{\varepsilon }(s) = \int_0^1 [ \, G_x ( X^{*} (s) + \lambda ( X_{\varepsilon }(s) - X^{*}(s) ) ) - G_x (X^{*}(s))  \, ]\, p(s) \, d\lambda ,$$
we get
\begin{eqnarray*}\label{eq7:lem3}
&& \hspace{-0.6cm} \xi_{\varepsilon } (t) = \xi_{\varepsilon } (t_0 + \varepsilon ) + \int_{t_0 + \varepsilon }^{t} \int_0^1 F_x ( X^{*} (s) + \lambda ( X_{\varepsilon }(s) - X^{*}(s) ) , u^{*} (s) ) \, \xi_{\varepsilon } (s)  \, d\lambda \, ds \nonumber \\
&& \hspace{0.75cm} + \, \int_{t_0 + \varepsilon }^{t} \int_0^1 G_x ( X^{*} (s) + \lambda ( X_{\varepsilon }(s) - X^{*}(s) ) ) \, \xi_{\varepsilon } (s)  \, d\lambda \, d M(s) \nonumber \\
&& + \, \int_{t_0 + \varepsilon }^{t} \int_0^1 [ \, F_x ( X^{*} (s) + \lambda ( X_{\varepsilon }(s) - X^{*}(s) ) , u^{*} (s) ) \nonumber \\
&& \hspace{2.5in} - F_x (X^{*}(s) , u^{*} (s) )  \, ]\, p(s) \, d\lambda \, ds \nonumber \\
&& \hspace{3in} + \, \int_{t_0 + \varepsilon }^{t} \tilde{\Phi}_{\varepsilon } (s) \, dM(s) ,
\end{eqnarray*}
for all $t \in [t_0 + \varepsilon , T] .$
Hence by making use of the isometry property (\ref{isometry property}) it holds $\forall \; t \in [t_0 + \varepsilon , T] ,$
\begin{eqnarray}\label{eq:for Gronwall-lem3}
&& \hspace{-0.75cm} \mathbb{E}\, [ \, |\, \xi_{\varepsilon } (t)  \, |^2 \, ] \leq 5\, \mathbb{E}\, [ \, |\, \xi_{\varepsilon } (t_0 + \varepsilon)  \, |^2 \, ] + 5 \, ( C_1 + C_2 ) \, \int_{t_0 + \varepsilon }^{t} \mathbb{E}\,  [ \, | \xi_{\varepsilon } (s)  \, |^2 \, ] \, ds
\nonumber \\
&&  \hspace{-0.75cm} + \, 5 \, \mathbb{E}\, \Big{[} \, \int_{t_0 }^T |\, \int_0^1 \big{(} \, F_x ( X^{*} (s) + \lambda ( X_{\varepsilon }(s) - X^{*}(s) ) , u^{*} (s) )
\nonumber \\
&& \hspace{2.5in} - F_x (X^{*}(s) , u^{*} (s) )  \, \big{)} p(s) \, d\lambda \, ds \, | \, \Big{]}^2
\nonumber \\
&&  \hspace{2in} + \, 5 \, \mathbb{E}\,  [ \, \int_{t_0}^T ||  \tilde{\Phi}_{\varepsilon } (s) \, \mathcal{Q}^{1/2} (s) ||^2_{2}  \, d s \, ] .
\end{eqnarray}
But as done for the second equality and first inequality in (\ref{eq:est1a-G-2}) we can derive easily that
\begin{eqnarray}\label{eq8:lem3}
&& \hspace{-0.5cm} \mathbb{E}\,  [ \, \int_{t_0}^T ||  \tilde{\Phi}_{\varepsilon } (s) \, \mathcal{Q}^{1/2} (s) ||^2_{2}  \, d s \, ] \nonumber \\
&& \hspace{-0.5cm} = \; \mathbb{E}\, [ \, \int_{t_0}^t || \tilde{\Phi}_{\varepsilon } (s) \mathcal{Q}^{1/2} (s) ||_2^2 \, ds \, ] \nonumber \\
&& \hspace{-0.5cm} \leq \; \mathbb{E}\, [ \, \int_{t_0}^t || \tilde{\Phi}_{\varepsilon } (s) ||_{L_{\mathcal{Q}} (K)}^2 \, ds \, ] \nonumber \\
&& \hspace{-0.5cm} = \; \mathbb{E}\, [ \, \int_{t_0}^t || \int_0^1 [ G_x ( X^{*} (s) + \lambda ( X_{\varepsilon }(s) - X^{*}(s) ) ) - G_x (X^{*}(s))  ] p(s) \, d\lambda ||_{L_{\mathcal{Q}} (K)}^2 ds ] \nonumber \\
&& \hspace{-0.5cm} \leq \; \mathbb{E}\, [ \, \int_{t_0}^t \int_0^1 || G_x ( X^{*} (s) + \lambda ( X_{\varepsilon }(s) - X^{*}(s) ) )
\nonumber \\
&& \hspace{2.25in}
- G_x (X^{*}(s))  \, ]\, p(s) ||_{L_{\mathcal{Q}} (K)}^2  \, d\lambda \, ds \, ] .
\end{eqnarray}
Therefore, from Lemma~\ref{lem:2nd estimate}, the continuity and boundedness of $G_x$ in (E1),
Lemma~\ref{lem:1st estimate} and the dominated convergence theorem to get that the last term in the right hand side of (\ref{eq8:lem3}) goes to $0$ as $\varepsilon \rightarrow 0 .$

Similarly, the third term in the right hand side of (\ref{eq:for Gronwall-lem3}) converges also to $0$ as $\varepsilon \rightarrow 0 .$

Finally, by applying Gronwall's inequality to (\ref{eq:for Gronwall-lem3}), and using (\ref{eq5:lem3})-(\ref{eq8:lem3}), we deduce that
\begin{equation*}\label{eq9:lem3}
\sup_{t_0 + \varepsilon \leq t \leq T } \, \mathbb{E}\, [ \, |\, \xi_{\varepsilon } (t)  \, |^2 \, ] \rightarrow 0 \;\; \text{as}\; \; \varepsilon \rightarrow 0 ,
\end{equation*}
which proves the lemma.
\end{proof}

\begin{lem}\label{lem:4th estimate}
Assume (E1) and (E2). Let $\zeta$ be the solution of the equation:
\begin{eqnarray*}\label{eq1:lem4}
 \left\{ \begin{array}{ll}
               d \zeta (t) = \ell_x (X^{*}(t) , u^{*} (t) ) p (t) dt , \hspace{0.5cm} t_0 < t \leq T,\\
               \, \zeta (t_0 ) = \ell (X^{*}(t_0) , v ) - \ell (X^{*}(t_0) , u^{*} (t_0) ) .
         \end{array}
 \right.
\end{eqnarray*}
Then
\begin{equation*}\label{eq2:lem4}
\lim_{\varepsilon \rightarrow 0 } \, \mathbb{E}\, \Big{[} \, \big{|} \, \frac{1}{\varepsilon } \int_{t_0 }^T \big{(} \ell ( X_{\varepsilon }(t) , u_{\varepsilon }(t) ) - \ell ( X^{*}(t) , u^{*} (t) ) \big{)} dt - \zeta (t) \, \big{|}^2 \, \Big{]} = 0 .
\end{equation*}
\end{lem}
\begin{proof}
Let $$\eta_{\varepsilon }(t) = \frac{1}{\varepsilon } \int_{t_0 }^t \big{(} \, \ell ( X_{\varepsilon }(t) , u_{\varepsilon }(t) ) - \ell ( X^{*}(t) , u^{*} (t) ) \big{)} dt - \zeta (T) ,$$ for $t \in [ t_0 , T ] .$ Then $\eta_{\varepsilon }(t_0) = - \big{(} \ell (X^{*}(t_0) , v ) - \ell (X^{*}(t_0) , u^{*} (t_0) ) \big{)} .$
So one can proceed easily as done in the proof of Lemma~\ref{lem:3rd estimate} to show that ${ \mathbb{E}\, [ \, | \, \eta_{\varepsilon }(T) \, |^2 \, ] \rightarrow 0 , }$ though this case is rather simpler.
\end{proof}

\bigskip

Let us now for a $C^1$ mapping $\Psi : K \rightarrow \mathbb{R}$ denote by $\nabla \Psi$ to the gradient of $\Psi ,$ which is
defined, by using the directional derivative $D\Psi (x) (k)$ of
$\Psi$ at a point $x \in K$ in the direction of $k \in K ,$ as
$\big{<} \nabla \Psi (x) , k \big{>} = D\Psi (x) (k) \, ( = \Psi_x (k) ).$ We shall sometimes write $\nabla_x \Psi $ for $\nabla \Psi (x) .$

\begin{cor}\label{cor1}
Under the assumptions of Lemma~\ref{lem:4th estimate}
\begin{equation}\label{eq1:cor1}
\frac{d}{d \varepsilon } \, J ( u_{\varepsilon } (\cdot ) ) \big{|}_{ \varepsilon = 0}  = \mathbb{E}\, [ \, \big{<} \nabla \, h ( X^{*}(T) ) , \, p (T) \big{>} + \zeta (T) \, ] .
\end{equation}
\end{cor}
\begin{proof}
Note that from the definition of the cost functional in (\ref{cost functional}) we see that
\begin{eqnarray*}\label{eq2:cor1}
&& \hspace{-0.5cm} \frac{1}{\varepsilon }~ \big{[} J ( u_{\varepsilon } (\cdot ) ) - J ( u^{*} (\cdot ) ) \big{]} = \frac{1}{\varepsilon }~ \mathbb{E}\, \Big{[} h ( X_{\varepsilon } (T) ) - h ( X^{*}(T) ) \nonumber \\ &&  \hspace{1.80in} + \, \int_{t_0 }^T \big{(} \ell (X_{\varepsilon } (s) , u_{\varepsilon } (s) ) -
\ell (X^{*} (s) , u^{*} (s) ) \big{)} ds \, \Big{]} \nonumber \\
&& = \mathbb{E}\, \Big{[} \, \int_0^1 h_x (X^{*} (T) + \lambda ( X_{\varepsilon }(T) - X^{*}(T) ) ) \, \frac{( X_{\varepsilon }(T) - X^{*}(T) )}{\varepsilon } \, d\lambda \nonumber \\ && \hspace{1.75in} + \; \frac{1}{\varepsilon }~  \int_{t_0 }^T \big{(} \ell (X_{\varepsilon } (s) , u_{\varepsilon } (s) ) - \ell (X^{*} (s) , u^{*} (s) ) \big{)} ds \, \Big{]}.
\end{eqnarray*}
Now let $\varepsilon \rightarrow 0 $ and use the properties of $h$ in (E1), Lemma~\ref{lem:3rd estimate} and Lemma~\ref{lem:4th estimate} to get (\ref{eq1:cor1}).
\end{proof}

\section{Maximum principle}\label{sec5}
The maximum principle is a good tool for studying the optimality of controlled SDEs like (\ref{forward-see}) since in fact the dynamic programming approach for similar optimal control problems require usually a Markov property to be satisfied by the solution of (\ref{forward-see}), cf. for instance \cite[Chapter~4]{[Y-Z]}. But this property does not hold in general especially when the driving noise is a martingale.

\smallskip

Let us recall the SDE~(\ref{forward-see}) and the mappings in (\ref{cost functional}), and define the \emph{Hamiltonian}
${ H:[ 0 , T ] \times \Omega \times K \times \mathcal{O} \times K \times L_2 (K) \rightarrow
\mathbb{R} }$ for $( t , \omega , x , u , y , z ) \in [ 0 , T ] \times \Omega \times K \times \mathcal{O}
\times K \times L_2 (K) $ by
\begin{equation}\label{defn of H}
 H ( t , \omega , x , u , y , z ) := \ell (x , u ) + \big<
F (x , u ) \, , y \big> + \big< G (x) \mathcal{Q}^{1/2} (t , \omega ) \, , z
\big>_2 \, .
\end{equation}

The adjoint equation of (\ref{forward-see}) is the following BSDE:
\begin{eqnarray}\label{adjoint-bse}
 \left\{ \begin{array}{ll}
             -\, d Y^{u (\cdot )}(t) = & \nabla_{x} H
             ( t , X^{u (\cdot )} (t), u (t), Y^{u (\cdot )} (t) , Z^{u (\cdot )}
              (t) \mathcal{Q}^{1/2} (t) ) \,
              dt  \\& ~ \hspace{1.5cm} - Z^{u (\cdot )} (t)\, d M(t) - d N^{u (\cdot )} (t)
             , \; \; \; t_0 \leq t < T , \\
             \; \; \; \, \, Y^{u (\cdot )} (T) = & \nabla h (X^{u (\cdot )} (T)) .
         \end{array}
 \right.
\end{eqnarray}

The following theorem gives the solution to BSDE~(\ref{adjoint-bse}) in the sense that there exists a triple
$( Y^{u (\cdot )} , Z^{u (\cdot )} , N^{u (\cdot )} )$ in ${L^2_{\mathcal{F}} ( 0 , T ; K )\times \Lambda^2 ( K ; \mathcal{P} ,
M ) \times \mathcal{M}^{2 , c}_{[ 0 , T] } (K) }$ such that the
following equality holds $a.s. $ for all $t \in [ 0 , T ] , \; N (0)
= 0$ and $N$ is VSO to $M$:
\begin{eqnarray*}
Y^{u (\cdot )} (t) &=& \xi + \int_t^T \nabla_{x} H
             ( s , X^{u (\cdot )} (s), u (s), Y^{u (\cdot )} (s) , Z^{u (\cdot )}
              (s) \mathcal{Q}^{1/2} (s) ) \,
              ds  \nonumber \\
&& - \int_t^T Z^{u (\cdot )} (s ) d M (s) - \int_t^T d N^{u (\cdot )} (s) .
\end{eqnarray*}
\begin{thm}\label{th:solution of adjointeqn}
Assume that (E1)--(E2) hold. Then there exists a unique solution
$( Y^{u (\cdot )} , Z^{u (\cdot )} , N^{u (\cdot )} )$ of the
BSDE~(\ref{adjoint-bse}).
\end{thm}

For the proof of this theorem one can see \cite{[Alh-Stoc09]}.

\bigskip

We shall denote briefly the solution of (\ref{adjoint-bse}), which corresponds to the optimal control $u^{*} ( \cdot )$ by $(Y^{*} , Z^{*} , N^{*}) .$

In the following lemma we shall try to compute $\mathbb{E}\, [ \, \big{<} Y^{*} (T) , \, p (T) \big{>}  \, ] .$
\begin{lem}\label{lem:formula of Y and p}
\begin{eqnarray}\label{Ito formula of Y and p}
&& \hspace{-1.5cm} \mathbb{E} \, [ \, \big< \, Y^{*} (T) , p(T) \, \big> \, ] = - \; \mathbb{E} \, \big{[} \, \int_{t_0 }^T \ell_x (X^{*} (s) , u^{*} (s) ) p(s) \, ds \, \big{]} \nonumber \\
&& \hspace{2cm} + \; \mathbb{E} \, \big{[} \, \big< Y^{*}(t_0 ) , F ( X^{*}(t_0 ) , v) - F ( X^{*}(t_0 ) , u^{*} (t_0 ) \big> \, \big{]} .
\end{eqnarray}
\end{lem}
\begin{proof}
Use It\^o's formula together to compute $d \, \big< Y^{*} (t) , p(t) \big> $ for $t \in [t_0 , T] ,$ and use the facts that
\begin{eqnarray*}
&& \hspace{-1.5cm}
\int_{t_0}^T \big< \, p(s) \, , \nabla_x  H ( s , X^{*} (s), u^{*} (s), Y^{*} (s) , Z^{*} (s) \mathcal{Q}^{1/2}(s) ) \big> \, ds
\\ &=&  \int_{t_0}^T \Big{[} \, \ell_x (X^{*} (s) , u^{*} (s) ) p(s) + \big< \, F_x (X^{*} (s) , u^{*} (s) ) p(s) \, , Y^{*}(s) \, \big> \Big{]}   \, ds
\\ && + \, \int_{t_0}^T \big{<} \, G_x ( X^{*} (s) ) p(s) \mathcal{Q}^{1/2} (s) \, , Z^{*}(s) \mathcal{Q}^{1/2} (s) \, \big{>}_2 \, ds ,
\end{eqnarray*}
which is easily seen from (\ref{defn of H}).
\end{proof}

\bigskip

Now we state our main result of this section.
\begin{thm}\label{main thm}
Suppose (E1)--(E2). If $( X^{*} , u^{*} (\cdot ) )$ is an optimal
pair for the problem (\ref{forward-see})-(\ref{value-function}), then there exists a unique solution $( Y^{*},
Z^{*} , N^{*} )$ to the corresponding BSDE~(\ref{adjoint-bse}) such that the following inequality holds:
\begin{eqnarray}\label{ineq1:main}
&& \hspace{-1cm} H ( t , X^{*} (t ) , v , Y^{*} (t ) , Z^{*} (t ) \mathcal{Q}^{1/2} (t ) )
 \nonumber \\ && \hspace{1.5cm} \geq \, H ( t , X^{*} (t ) , u^{*} (t ) , Y^{*} (t) , Z^{*} (t ) \mathcal{Q}^{1/2} (t ) )  \\
&& \hspace{2.5in}  \text{a.e.} \; t \in [0 , T],\; \text{a.s.} \; \forall \; v \in U .  \nonumber
\end{eqnarray}
\end{thm}
\begin{proof}
We note that since $u^{*} (\cdot )$ is optimal, $\frac{d}{d \varepsilon } \, J ( u_{\varepsilon } (\cdot ) )|_{\varepsilon = 0 } \, \geq 0 ,$ which implies by using Corollary~\ref{cor1} that
\begin{equation}\label{eq:useful-ineq}
\mathbb{E}\, [ \, \big{<} Y^{*} (T) , \, p (T) \big{>} + \zeta (T) \, ] \geq 0 .
\end{equation}

On other hand by applying (\ref{eq:useful-ineq}) and Lemma~\ref{lem:formula of Y and p} one sees that
\begin{eqnarray}\label{eq1:mainthm}
&& \hspace{-1cm} 0 \leq - \; \mathbb{E} \, [ \, \int_{t_0 }^T \ell_x (X^{*} (s) , u^{*} (s) ) p(s) \, ds \, ]  \nonumber \\
&& \hspace{1cm} + \;   \mathbb{E} \, [ \, \big< Y^{*}(t_0 ) , F ( X^{*}(t_0 ) , v) - F ( X^{*}(t_0 ) , u^{*} (t_0 ) \big> + \zeta (T) \, ] .
\end{eqnarray}
But
$$ \zeta (T) = \zeta (t_0 ) + \int_{t_0 }^T \ell_x (X^{*} (s) , u^{*} (s) ) p(s) \, ds $$
and
\begin{eqnarray*}
&& \hspace{-1.5cm} H ( t_0 , X^{*} (t_0 ) , v , Y^{*} (t_0 ) , Z^{*} (t_0 ) \mathcal{Q}^{1/2} (t_0 ) )
\\ && \hspace{1cm} - H ( t_0 , X^{*} (t_0 ) , u^{*} (t_0 ) , Y^{*} (t_0) , Z^{*} (t_0 ) \mathcal{Q}^{1/2} (t_0 ) )
\\ &&  = \zeta (t_0 ) + \big< Y^{*}(t_0 ) , F ( X^{*}(t_0 ) , v ) - F ( X^{*}(t_0 ) , u^{*} (t_0 ) ) \big> .
\end{eqnarray*}
Hence (\ref{eq1:mainthm}) becomes
\begin{eqnarray}\label{eq2:mainthm}
&& 0 \leq \mathbb{E} \; [ \, H ( t_0 , X^{*} (t_0 ) , v , Y^{*} (t_0 ) , Z^{*} (t_0 ) \mathcal{Q}^{1/2} (t_0 ) )
 \nonumber \\ && \hspace{1.5cm} - \, H ( t_0 , X^{*} (t_0 ) , u^{*} (t_0 ) , Y^{*} (t_0) , Z^{*} (t_0 ) \mathcal{Q}^{1/2} (t_0 ) ) \, ] .
\end{eqnarray}

Now varying $t_0$ as in (\ref{eq:t0}) shows that (\ref{eq2:mainthm}) holds for $a.e. \; t . ,$ and so by arguing for instance as in \cite[P. 19]{[Be82]} we obtain easily (\ref{ineq1:main}).
\end{proof}
%
%

\bigskip

\begin{rk}
Let us assume for example that the space $K$ in Theorem~\ref{th:solution of adjointeqn} is the real space $\mathbb{R}$ and $M$ is the martingale
given by the formula $$M(t) = \int_0^t \alpha (s) dB(s) , \hspace{1cm} t \in [0 , T] ,$$ for some $\alpha \in L^2_{\mathcal{F}} (0 , T ; \mathbb{R} ) $ and a one dimensional Brownian motion $B$. If
$\alpha(s) > 0$ for each $s ,$ then $\mathcal{F}_t (M) = \mathcal{F}_t (B)$ for each $t ,$
where $$\mathcal{F}_t (R) = \sigma \{ R(s) , 0 \leq s \leq t \} $$ for $R = M , B .$
Consequently, by applying the unique representation property for martingales with respect to $\{ \mathcal{F}_t (M) , \; t \geq 0 \}$ or larger filtration in \cite[Theorem 2.2]{[Alh-Stoc09]} or \cite{[Alh-ROSE10]} and
the Brownian martingale representation theorem as e.g. in [14, Theorem 3.4, P. 200], we deduce that the martingale $N^{ u ( \cdot ) } $ in (\ref{adjoint-bse}) vanishes almost surely if the filtration furnished for the SDE~(\ref{forward-see}) is $\{ \mathcal{F}_t (M) , \; 0 \leq t \leq T \} .$ This result follows from the construction of the solution of the BSDE~(\ref{adjoint-bse}). More details on this matter can be found in \cite[Section~3]{[Alh-Stoc09]}. As a result, in this particular case BSDE~(\ref{adjoint-bse}) fits well with those BSDEs studied by Pardoux \& Peng in \cite{[Pa-Pe90]}, but with the variable $\alpha Z$ replacing $Z$ there.

Thus in particular we conclude that many of the applications  of BSDEs, which were studied in the literature, to both stochastic optimal control and finance (e.g. \cite{[Zhou93]} and the references therein) can be applied directly or after slight modification to work here for BSDEs driven by martingales. For example we refer the reader to \cite{[Marie09]} for financial application. Another interesting case can be found in \cite{[Bal-Pard05]}.

On the other hand, in this respect we shall present an example (see Example~\ref{Example2}) in Section~\ref{sec6}, by modifying an interesting example due to Bensoussan \cite{[Be82]}.
\end{rk}

\section{Sufficient conditions for optimality}\label{sec6}
In the previous two sections we derived Pontyagin's maximum principle which gives necessary conditions for optimality for the control problem (\ref{forward-see})-(\ref{value-function}). In the following theorem, if we have also a convexity assumption on the control domain $U ,$ we shall obtain sufficient conditions for optimality of this optimal control problem. This concerned result is a variation of Theorem~4.2 in \cite{[Alh-Int10]}.
\begin{thm}\label{mainthm2}
Assume (E1) and, for a given $u^* (\cdot ) \in \mathcal{U}_{ad} ,$ let $X^{*}$ and $( Y^{*},
Z^{*} , N^{*} )$ be the corresponding solutions of
equations (\ref{forward-see}) and (\ref{adjoint-bse}) respectively.
Suppose that the following conditions hold:\\
$(i)$ $U$ is a convex domain in $\mathcal{O} ,$ $h$ is convex, \\
$(ii)$ $(x , v) \mapsto H ( t , x , v , Y^{*}(t) , Z^{*}(t) \mathcal{Q}^{1/2} (t) )$ is convex for all $t \in [0 , T]$ \,a.s.,
\vspace{-0.35cm}
\noindent\begin{eqnarray*} & & \hspace{-0.70in} (iii)\; H ( t , X^{*} (t) , u^{*} (t) ,
Y^{*} (t) , Z^{*} (t) \mathcal{Q}^{1/2} (t) )
\\
& & \hspace{1.5in} = \min_{v \in U } \; H ( t ,
X^{*}(t) , v , Y^{*}(t) ,
Z^{*}(t) \mathcal{Q}^{1/2} (t) )
\end{eqnarray*}
for a.e. $t \in [0 , T]$ a.s.

Then $( X^{*} , u^{*} (\cdot ) )$ is an optimal
pair for the control problem (\ref{forward-see})-(\ref{value-function}).
\end{thm}
\begin{proof}
Let $u (\cdot ) \in \mathcal{U}_{ad} .$ Consider the following definitions:
\[ I_1 : = \mathbb{E} \; \big[ \; \int_0^T
\big( \ell ( X^{*}(t) , u^{*} (t ) ) - \ell(X^{u
(\cdot )}(t) , u (t ) ) \big) dt \; \big] \]
and
\[ I_2 : = \mathbb{E} \; [ \, h (  X^{*}(T) )  - h (
X^{u (\cdot )}(T) ) \, ] . \]
Then readily
\begin{eqnarray}\label{eq1:mainthm2}
 J ( u^{*} (\cdot ) ) - J ( u  (\cdot ) ) = I_1 + I_2 .
\end{eqnarray}

Let us define
\begin{eqnarray*}
&& I_3 : = \mathbb{E} \; \big[\; \int_0^T
 \big{(} H ( t , X^{*} (t) , u^{*} (t) , Y^{*} (t) , Z^{*} (t) \, \mathcal{Q}^{1/2} (t) )
\\ && \hspace{1.3in} - \, H ( t , X^{u (\cdot ) } (t) , u (t) , Y^{*} (t) , Z^{*} (t) \, \mathcal{Q}^{1/2} (t) ) \big{)}  dt \; \big] ,
\end{eqnarray*}
\begin{eqnarray*}
I_4 := \mathbb{E} \; \big[ \, \int_0^T \big{<} F ( X^{*}
(t) , u^{*} (t) ) - F ( X^{u (\cdot ) } (t) , u (t) ) \, ,
Y^{*} (t) \big{>} \, dt \, \big] ,
\end{eqnarray*}
\begin{eqnarray*}
I_5 := \mathbb{E} \; \big[\, \int_0^T \big{<} \big( G ( X^{u^{*}
(\cdot ) } (t) ) - G ( X^{u (\cdot ) } (t) ) \big)
\mathcal{Q}^{1/2} (t) \, , Z^{*}(t)
\mathcal{Q}^{1/2} (t) \big{>}_2 \, dt \, \big] ,
\end{eqnarray*}
and
\begin{eqnarray*}
&& \hspace{-1.2cm} I_6 := \mathbb{E} \; \Big[ \; \int_0^T \big{<} \nabla_{x} H ( t ,
X^{*} (t) , u^{*} (t) , Y^{*} (t)
, Z^{u^{*} (\cdot )} (t) \mathcal{Q}^{1/2} (t) ) \, , \\ && \hspace{3in} X^{*}(t) - X^{u
(\cdot ) }(t) \big{>} \, dt \; \Big] .
\end{eqnarray*}

From the definition of $H$ in (\ref{defn of H}) we get
\begin{equation}\label{eq2:mainthm2}
 I_1 = I_3 - I_4 - I_5 .
\end{equation}

On the other hand, from the convexity of $h$ in condition (ii) it follows
\begin{eqnarray*}
h (  X^{*}(T) )  - h ( X^{u (\cdot )}(T) )
\leq \; \big{<} \, \nabla h (  X^{*}(T) ) \, ,
X^{*}(T) - X^{u (\cdot )}(T)
 \, \big{>} \; \; a.s.,
\end{eqnarray*}
which implies that
\begin{equation}\label{eq3:mainthm2}
 I_2 \leq \mathbb{E} \; [ \, \big{<} \; Y^{*}(T) \; ,
X^{*} (T) - X^{u }(T) \; \big{>} \, ] .
\end{equation}

Next by applying It\^{o}'s formula to compute $d \, \big{<} Y^{*}
(t) \, , X^{*}(t) - X^{u (\cdot ) }(t) \big{>} $ and using equations (\ref{adjoint-bse}) and (\ref{forward-see}) we find with the help of (\ref{eq3:mainthm2}) that
\begin{equation}\label{eq4:mainthm2}
I_2 \leq I_4 + I_5 - I_6 \, .
\end{equation}
Consequently, by considering (\ref{eq1:mainthm2}), (\ref{eq2:mainthm2}) and (\ref{eq4:mainthm2}) it follows that
\begin{equation}\label{eq5:mainthm2}
 J ( u^{*} (\cdot ) ) - J ( u  (\cdot ) ) \leq I_3 - I_6 .
\end{equation}

On the other hand, from the convexity property of the mapping
${ (x , v) \mapsto H ( t , x , u , Y^{*}(t) , Z^{*}(t)  \mathcal{Q}^{1/2} (t) ) }$ in assumption (iii) the following inequality holds a.s.:
\begin{eqnarray*}\label{eq6:mainthm2}
& & \hspace{-0.50cm} \int_0^T \Big{(} H ( t , X^{*} (t) , u^{*}
(t) , Y^{*} (t) , Z^{*} (t)  \mathcal{Q}^{1/2} (t) ) \\ & &
\hspace{1.5in} - \; H ( t , X^{u (\cdot ) } (t) , u (t) , Y^{*}
(t) , Z^{*} (t)  \mathcal{Q}^{1/2} (t) ) \Big{)} \; dt  \\ & & \hspace{-0.5cm} \leq \,
\int_0^T \big{<} \; \nabla_{x} H ( t , X^{*}(t),
u^{*} (t), Y^{*} (t) , Z^{*} (t)  \mathcal{Q}^{1/2} (t) )
\, , \, \\ && \hspace{3.2in} X^{*} (t) - X^{u (\cdot ) } (t) \; \big{>}
\; dt \\ & & \hspace{-0.5cm} + \int_0^T \big{<} \; \nabla_{u } H ( t ,
X^{*} (t), u^{*} (t), Y^{*} (t) ,
Z^{*}(t)  \mathcal{Q}^{1/2} (t) ) \, , \, u^{*} (t) - u (t) \, \big{>}_{\mathcal{O}}
\, dt .
\end{eqnarray*}
As a result
\begin{equation}\label{eq7:mainthm2}
I_3 \leq I_6 + I_7 ,
\end{equation}
where
\begin{eqnarray*}
&& I_7 = \mathbb{E} \, \Big[ \; \int_0^T \big{<} \; \nabla_{u } H ( t ,
X^{*} (t), u^{*} (t), Y^{*} (t) ,
Z^{*}(t) \mathcal{Q}^{1/2} (t) ) \, , \\ && \hspace{3.6in} u^{*} (t) - u (t) \;
\big{>}_{\mathcal{O} } \, dt \; \Big] .
\end{eqnarray*}
Since $v \mapsto H ( t , X^{*} (t), v, Y^{*} (t) , Z^{*}(t) \mathcal{Q}^{1/2} (t) )$ is minimum at $v = u^{*} (t)$ by the minimum condition (iii), we have
$$\big{<} \; \nabla_{u } H ( t , X^{*} (t), u^{*} (t), Y^{*} (t) , Z^{*}(t) \mathcal{Q}^{1/2} (t) ) \, , u^{*} (t) - u (t) \;
\big{>}_{\mathcal{O} } \leq 0 .$$
Therefore $I_7 \leq 0 ,$ which by (\ref{eq7:mainthm2}) implies that
$I_3 - I_6 \leq 0 .$ So (\ref{eq5:mainthm2}) becomes $$ J ( u^{*} (\cdot ) ) - J ( u  (\cdot ) ) \leq 0 .$$

Now since $u (\cdot ) \in \mathcal{U}_{ad} $ is arbitrary, this inequality proves that $( X^{*} , u^{*} (\cdot ) )$ is an optimal
pair for the control problem (\ref{forward-see})-(\ref{value-function}) as required.
\end{proof}

\begin{ex}\label{Example1}
Let $m$ be a continuous square integrable one dimensional martingale with respect to $\{ \mathcal{F}_t \}_t$ such that
$ <m>_t \, = \int_0^t \alpha(s) ds$ $\forall ~ 0 \leq t \leq T$ for some continuous $\alpha: [0 , T] \rightarrow (0 , \infty ) .$
Consider $M(t) = \beta \, m(t) ( = \int_0^t \beta \, d m(s) ), $ with $\beta \neq 0$ being a fixed
element of $K .$ Then $M \in \mathcal{M}^{2 , c} (K) $ and $<<M>>_t $ equals
$\widetilde{\beta\otimes \beta } \;\int_0^t \alpha(s) ds , $ where
$\widetilde{\beta\otimes \beta }$ is the identification of
$\beta\otimes \beta$ in $L_1 (K) , $ that is $(\widetilde{\beta\otimes
\beta }) (k) = \big{<} \beta , k \big{>} \, \beta , \; k \in K .$ Also $<M>_t \; = | \beta |^2 \int_0^t \alpha (s) \, ds .$
Now letting $\mathcal{Q} (t) = \widetilde{ \beta\otimes \beta } \; \alpha(t) $ yields
that $<<M>>_t \; = \int_0^t \mathcal{Q} (s) \, ds .$ This process $\mathcal{Q}(\cdot )$ is bounded since
$\mathcal{Q} (t) \leq \mathcal{Q} \; \; \forall \; t ,$
where $\mathcal{Q} = \widetilde{ \beta\otimes \beta } \; \displaystyle{\max_{0 \leq t \leq T } }\alpha(t) .$
It is also easy to see that $\mathcal{Q}^{1/2} (t) (k) = \frac{\big{<} \beta , k \big{>}
\, \beta}{| \beta |} \; \alpha^{1/2} (t) .$ In particular $\beta \in \mathcal{Q}^{1/2} (t) (K) .$

\smallskip

Let $K = L^2 ( \mathbb{R}^n ) .$ Let $M$ be the above martingale.
Suppose that $\mathcal{O} = K .$
Assume that $\tilde{G} \in L_{\mathcal{Q}} (K) $ or
even a bounded linear operator from $K$ into itself, and $\tilde{F} $ is a bounded
linear operator from $\mathcal{O}$ into $K .$ Let us consider the SDE:
\begin{eqnarray*}
 \left\{ \begin{array}{ll}
              d X(t) = \tilde{F} \, u (t) \; d t
              + \, \big{<} X (t ) \, , \beta \big{>} \;
              \tilde{G} \; d M(t) , \;\; t \in [0 , T] ,\\
             \; X(0) = x_0 \in K .
         \end{array}
 \right.
\end{eqnarray*}

For a given fixed element $c$ of $K$ we assume that the cost functional is given by the formula:
\begin{equation*}
J ( u (\cdot ) ) = \mathbb{E} \, [ \, \int_0^T  | u
(t) |^2 \, dt \, ] + \mathbb{E} \; [ \; \big{<} c \, , X (T)
\big{>} \, ] ,
\end{equation*}
and the value function is \[ J^{*} = \inf \{ J( u (\cdot
) ) : \; u (\cdot ) \in \mathcal{\mathcal{U}}_{ad} \} . \]

This control problem can be related to the control problem
(\ref{forward-see})-(\ref{value-function}) as follows.
We define
 \[ F ( x , u )
= \tilde{F} \, u , \; G (x) = \big{<} x \, , \beta \big{>} \,
\tilde{G} , \; \; \ell (x , u) = |u|^2 , \; \; \text{and} \; \; h (x) = \big{<} c \, , x \big{>} \, ,
\] where $(x , u ) \in K \times \mathcal{O} .$

The Hamiltonian then becomes the mapping
\[ H: [0 , T] \times \Omega \times K \times \mathcal{O} \times K \times L_2 (K) \rightarrow
\mathbb{R} , \]
\begin{equation*}
 H (t, x , u , y , z) = | u |^2 + \big<
\tilde{F} \, u  \, , y \big> + \big{<} x \, , \beta \big{>} \; \big<
\tilde{G} \, \mathcal{Q}^{1/2} (t) \, , z \big>_2 \, ,
\end{equation*}
$( t , x , u , y , z ) \in K \times \mathcal{O} \times K \times L_2 (K) .$

It is obvious that $H ( \cdot , \cdot , y , z )$ is convex with
respect to $(x , u ) $ for each $y$ and $z $ and $ \nabla_{x }
H ( t , x , u , y , z ) = \big< \tilde{G} \, \mathcal{Q}^{1/2} (t) \, ,
z \big> \; \beta .$

Next we consider the adjoint BSDE:
\begin{eqnarray*}
 \left\{ \begin{array}{ll}
             -\, d Y (t) = [\, \big< \tilde{G} \, \mathcal{Q}^{1/2} (t)
             \, , Z (t) \big>_2 \; \beta
             \; ]\; dt  - Z (t)\; d M(t) - d N (t), \\
             \; \; \; Y (T) = c .
          \end{array}
 \right.
\end{eqnarray*}
This BSDE attains an explicit solution $Y(t) = c
\, ,$ since $c$ is non-random. But this implies that $Z(t) = 0 $ and $N(t) = 0 $ for each
$t \in [ 0 , T ] .$

On the other hand, we note that the function $ \mathcal{O} \ni u \mapsto
H ( t , x , u , y , z ) \in \mathbb{R} $ attains its minimum
at $u = - \, \frac{1}{2} \, \tilde{F}^{*} \, y ,$ for fixed $(x ,
y , z) .$ So we choose our candidate for an optimal control as
\[ u^* (t , \omega ) = - \, \frac{1}{2} \,
\tilde{F}^{*} \, Y(t , \omega ) = - \, \frac{1}{2} \, \tilde{F}^{*} \,
c \; \; ( \in U := \mathcal{O} ) , \] .

With this choice all the requirements in
Theorem~\ref{mainthm2} are verified. Consequently
$u^* (\cdot )$ is an optimal control of this control problem with
an optimal solution $\hat{X} $ given by the solution of the
following closed loop equation:
\begin{eqnarray*}
 \left\{ \begin{array}{ll}
              d \hat{X}(t) = - \, \frac{1}{2} \, \tilde{F} \, \tilde{F}^{*} \, Y(t) \, d t
              + \, \big{<} \hat{X} (t ) \, , \beta \big{>} \;
              \tilde{G} \, d M(t) ,\\
             \; \hat{X}(0) = x_0 \in K .
         \end{array}
 \right.
\end{eqnarray*}

The value function takes the following value:
\begin{equation*}
J^{*} = \, \frac{1}{4} \;  | \tilde{F}^{*} c |^2 \, T + \mathbb{E} \;
[ \; \big{<} c \; , \hat{X} (T) \big{>} \; ] .
\end{equation*}

\end{ex}

\begin{rk}
It would be possible if we take $h (x) = |x|^2 , \; x \in K ,$ in the preceding example and proceeds as above. However if a result of existence and uniqueness os solutions to what we may call ``forward-backward stochastic differential
equations with martingale noise" holds, it should certainly be very useful
to deal with both this particular case and similar problems.
\end{rk}

\begin{ex}\label{Example2}
Let $\mathcal{O} = K .$ We are interested in the following linear quadratic example, which is gleaned from Bensoussan \cite[P. 33]{[Be82]}. Namely, we consider the SDE:
\begin{eqnarray}\label{ex2:forward-see}
 \left\{ \begin{array}{ll}
               d X(t) = ( A(t) X(t) + C(t) u (t) + f(t) )\, d t
              + \, ( B(t) X(t) + D(t) )\, d M(t) , \\
             \; X(0) = x_0 ,
         \end{array}
 \right.
\end{eqnarray}
where $B (t) x = \big< \gamma (t) \, , x \big> \, \tilde{G} (t) $ and $A , \gamma , C : [0,T] \times K \rightarrow K, \; f : [0,T] \rightarrow K , \; \tilde{G} , D: [0,T] \rightarrow L_{\mathcal{Q}}(K)$ are measurable and bounded mappings.

Let $P, Q : [0,T] \times K \rightarrow K, \; P_1: K \rightarrow K$ be measurable and bounded mappings. Assume that $P, P_1$ are symmetric non-negative definite, and $Q$ is a symmetric positive definite and $Q^{-1}(t)$ is bounded.
For SDE~(\ref{ex2:forward-see}) we shall assume that the cost functional is
\begin{eqnarray}\label{ex2:cost functional}
&& \hspace{-1.5cm} J(u (\cdot ) ) = \mathbb{E} \, \Big[ \, \int_0^T  \big( \, \frac{1}{2} \,  \big< P(t) X^{u (\cdot ) } (t)  \, , X^{u (\cdot ) } (t)  \big>  + \frac{1}{2} \,  \big< Q(t) u (t) \, , u (t) \big> \, \big) \, dt \nonumber \\ &&  \hspace{2.25in} + \, \frac{1}{2} \,  \big< P_1 X^{u (\cdot ) } (T) \, , X^{u (\cdot ) } (T) \big> \, \Big] ,
\end{eqnarray}
for $u (\cdot ) \in \mathcal{U}_{ad} .$

\smallskip

The control problem now is to minimize (\ref{ex2:cost functional}) over the set $\mathcal{U}_{ad} $ and get an optimal control $u^{*} ( \cdot ) \in \mathcal{U}_{ad} ,$ that is
\begin{equation}\label{ex2:value-function}
J(u^{*} ( \cdot ) ) = \inf \{ J(u (\cdot ) ): \; u (\cdot ) \in \,
\mathcal{U}_{ad} \} .
\end{equation}

By recalling Remark~\ref{rk:dependence on t} we can consider this control problem (\ref{ex2:forward-see})-(\ref{ex2:value-function}) as a control problem of the type (\ref{forward-see})-(\ref{value-function}). To this end, we let $$F(t,x,u) = A(t) x + C(t) u + f(t) , $$
$$G(t, x) = \big< \gamma(t) \, , x \big> \, \tilde{G} (t) + D(t) ,$$
$$\ell (t,x,u) = \frac{1}{2} \,  \big< P(t) x \, , x \big> + \frac{1}{2} \,  \big< Q(t) u \, , u \big> ,$$
$$h(x) = \frac{1}{2} \,  \big< P_1 x \, , x \big> .$$
Then the Hamiltonian ${ H:[ 0 , T ] \times \Omega \times K \times K \times K \times L_2 (K) \rightarrow\mathbb{R} }$ is given by
\begin{eqnarray*}\label{ex2:defn of H}
&& \hspace{-0.5cm} H ( t , x , u , y , z ) = \ell (t, x , u ) + \big<
F (t, x , u ) \, , y \big> + \big< G (t,x) \mathcal{Q}^{1/2} (t) \, , z
\big>_2 \, \nonumber \\ && \hspace{2.25cm} = \, \frac{1}{2} \,  \big< P(t) x \, , x \big> + \frac{1}{2} \,  \big< Q(t) u \, , u \big> \nonumber \\ && \hspace{2.1in} + \,
\big< A(t) x + C(t) u + f(t) \, , y \big>
\nonumber \\ && \hspace{2.1in}
+ \, \big< \, ( \big< \gamma(t) \, , x \big> \, \tilde{G} (t) + D(t) ) \mathcal{Q}^{1/2} (t) \, , z
\big>_2 \, .
\end{eqnarray*}
We can compute $\nabla_x H$ directly to find that
$$\nabla_x H ( t , x , u , y , z ) = P(t) u + A^* (t) x + \big< \tilde{G} (t) \mathcal{Q}^{1/2} (t) \, , z \big>_2 \, \gamma(t) .$$
Hence the adjoint equation of (\ref{ex2:forward-see}) takes the following shape:
\begin{eqnarray*}\label{ex2:adjoint-bse}
 \left\{ \begin{array}{ll}
              -\, d Y^{u (\cdot )}(t) = & \Big( \, A^* (t) Y^{u (\cdot )}(t) + P(t) X^{u (\cdot )}(t)
             \\ ~ & \hspace{2cm}  + \big< \tilde{G} (t) \mathcal{Q}^{1/2} (t) \, , Z^{u (\cdot )} (t) \mathcal{Q}^{1/2} (t) \big>_2 \, \gamma(t) \, \Big) \, dt
             \\ ~ & \hspace{2in}  - Z^{u (\cdot )} (t) d M(t) - d N^{u (\cdot )} (t) , \\
             \; \; \; \, \, Y^{u (\cdot )} (T) = & P_1 X^{u (\cdot )} (T) .
         \end{array}
 \right.
\end{eqnarray*}

Now the maximum principle theorems (Theorem~\ref{main thm}, Theorem~\ref{mainthm2}) in this case hold readily if we consider Remark~\ref{rk:dependence on t} again, and yield eventually
$$C^* (t) Y^* (t) + Q(t) u^* (t) = 0 .$$
\end{ex}

\bigskip

\bigskip

{\bf Acknowledgement.}
This author would like to thank the Mathematics Institute, Warwick University, where part of this work was done,
for hospitality during the summer of 2011.

\fussy

\end{document}